\documentclass[a4paper, 12pt]{article}
\usepackage{amsfonts}

 \usepackage{mathrsfs,amsfonts,amsmath}
 \usepackage{color}

\makeatletter

\@addtoreset{equation}{section}
\makeatother
\newtheorem{Theorem}{Theorem}[section]
\newtheorem{Remark}[Theorem]{Remark}
\newtheorem{Lemma}[Theorem]{Lemma}
\newtheorem{Corollary}[Theorem]{Corollary}

\begin{document}

\title{\textbf{New sufficient conditions of existence, moment estimations and non confluence
for SDEs with non-Lipschitzian coefficients}}
\author{Guangqiang Lan\footnote{Corresponding author. Supported by China Scholarship Council,
National Natural Science Foundation of China (NSFC11026142) and Beijing Higher Education Young Elite Teacher Project (YETP0516). }
\\ \small School of Science, Beijing University of Chemical Technology, Beijing 100029, China
\\ \small Email: langq@mail.buct.edu.cn
\\  Jiang-Lun Wu
\\ \small Department of Mathematics, College of Science, Swansea University, Swansea SA2 8PP, UK
\\ \small Email: j.l.wu@swansea.ac.uk}
\date{}

\maketitle

\begin{abstract}
The object of the present paper is to find new sufficient conditions for the existence of unique strong solutions to
a class of (time-inhomogeneous) stochastic differential equations with random, non-Lipschitzian coefficients.
We give an example to show that our conditions are indeed weaker than those relevant conditions existing in the literature. We also derive moment estimations for the maximum process of the solution. Finally, we present a
sufficient condition to ensure the non confluence property of the solution of time-homogeneous SDE which, in
one dimension, is nothing but stochastic monotone property of the solution.
 \end{abstract}



\section{ Introduction and Main Results}

\noindent
The theory of stochastic differential equations (SDEs) has been very well developed since the seminal work
of the great mathematician Kiyosi It\^{o} in the mid 1940s. Fundamental conditions like linear growth and Lipschitzian
type conditions on the both drift and diffusion coefficients to ensure the existence and uniqueness of solutions of
SDEs with any given initial data. The proofs are either based on Picard iteration (see. e.g., \cite{Ikeda}) or via
martingale problem formulation (cf. \cite{Stroock}). Since the remarkable paper \cite{Fang}, SDEs (as well as
stochastic functional differential equations) with non-Lipschitzian coefficients have received
much attention widely, see, e.g., \cite{Lan,Renesse,Hof,Wang}, just mention a few. In the present paper, we aim to
issue new sufficient conditions for the existence and uniqueness of strong solutions to SDEs with random non-
Lipschitzian coefficients. We also derive the moment estimation of the solution. Furthermore, we will give sufficient conditions for the non confluence property (also known as non contact property, cf. \cite{Yamada})  of the solution
of time-homogeneous SDE.

Given a probability space $(\Omega,\mathscr{F},P)$ endowed with a complete filtration $(\mathscr{F}_t)_{t\geq 0}$.
Let $d,m\in\mathbb{N}$ be arbitrarily fixed. We are concerned with the following stochastic differential equation

\begin{equation}\label{sdelan}dX_t(\omega)=\sigma(t,\omega,X_t)dB_t(\omega)+b(t,\omega,X_t)dt\
X_0(\omega)=x_0,a.s. \end{equation}
where the initial $x_0\in \mathbb{R}^d, (B_t)_{t\geq0}$ is an $m$-dimensional standard $\mathscr{F}_t$-Brownian motion,
and $\sigma:(t,\omega,x)\in[0,\infty)\times\Omega\times\mathbb{R}^d\mapsto\sigma(t,\omega,x)\in\mathbb{R}^d\otimes\mathbb{R}^m$
and $b:(t,\omega,x)\in[0,\infty)\times\Omega\times\mathbb{R}^d\mapsto b(t,\omega,x)\in\mathbb{R}^d$ are progressively
measurable, respectively, continuous with respect to the third variable $x$.

In order that the integrals in the definition of the solutions of the equation (\ref{sdelan}) are well-defined, we make
the following assumption which is enforced throughout the paper

\begin{equation}\label{t1}\mathbb{E}\int_0^T\sup_{|x|\le R}(|b(s,\cdot,x)|+||\sigma(s,\cdot,x)||^2)ds<\infty,\quad \forall T,R>0
\end{equation} where the norm $||\cdot||$ stands for the Hilbert-Schmidt norm $||\sigma||^2:=\sum\limits^d_{i=1}\sum\limits^m_{j=1}\sigma^2_{ij}$
for any $d\times m$-matrix $\sigma=(\sigma_{ij})\in\mathbb{R}^d\otimes\mathbb{R}^m$ and $|\cdot|$ denotes the usual Euclidean norm on $\mathbb{R}^d$.

Let us first discuss the sufficient conditions under which there is a unique strong solution of equation $(\ref{sdelan})$.
Fix $R>0$ arbitrarily, let $\eta_R:[0,1)\rightarrow\mathbb{R}_+$ be an increasing
continuous function on a very small interval $[0,\varepsilon_0)(0<\varepsilon_0\le c_0)$ which  satisfies

$$\eta_R(0)=0,\int_{0+}\frac{dx}{\eta_R(x)}=\infty.$$

Our first main result is the following

\begin{Theorem}\label{dingli1} Let $R>0$ be fixed arbitrarily.  Assume that for all $t\ge 0, \omega\in \Omega,$
$|x|\vee|y|\le R,$ the following locally weak monotonicity condition
\begin{equation}\label{t3}
||\sigma(t,\omega,x)-\sigma(t,\omega,y)||^2+2\langle x-y,b(t,\omega,x)-b(t,\omega,y)\rangle
\leq g(t,\omega)\eta_R(|x-y|^2)
\end{equation}
holds for $|x-y|\le c_0<1$, with $g$ being a progressively measurable and non negative function such that
$$\mathbb{E}\int_0^t g(s,\cdot)ds<\infty,\ t\geq 0.$$
Then, there is a unique strong solution of SDE $(\ref{sdelan})$.\end{Theorem}

We will use Euler's approximation method to prove Theorem
\ref{dingli1}. To this end, we give briefly the Euler's approximation for our SDE $(\ref{sdelan})$.
For any fixed $n,$ set $X^{(n)}(0):=x_0,$ for nonnegative integer $k$, and $t\in[\frac{k}{n},\frac{k+1}{n}) $,
define

\begin{equation}\label{t8} X^{(n)}(t):=X^{(n)}(\frac{k}{n})+\int_{\frac{k}{n}}^tb(s,\cdot,X^{(n)}(\frac{k}{n}))ds
+\int_{\frac{k}{n}}^t\sigma(s,\cdot,X^{(n)}(\frac{k}{n}))dB_s.
\end{equation}

To prove \textbf{Theorem \ref{dingli1}}, we need the following lemmas.

\begin{Lemma}\label{l1}

Assume that $(\ref{t3})$ holds. Then there is at least one uniformly convergent subsequence of $ X^{(n)}(\cdot)$.

\end{Lemma}

\begin{Lemma}\label{l2}

Assume that $(\ref{t3})$ holds. If $X(\cdot)$ is the limit process of a subsequence of $X^{(n)}(\cdot)$, then $X^{(n)}(\cdot)$ is a solution of equation $(\ref{sdelan})$.

\end{Lemma}

Our next main result concerns the non explosion property of the solution of SDE $(\ref{sdelan})$.  Let
$\gamma:[0,\infty)\rightarrow\mathbb{R}_+$ be a continuous, increasing function satisfying
$$(i)\
\lim_{x\rightarrow\infty}\gamma(x)=\infty;\ (ii)
\int_K^{\infty}\frac{dx}{\gamma(x)+1}=\infty.
$$

Our second main result is the following

\begin{Theorem}\label{dingli2}  Assume that there is a constant $K>0$ such that

\begin{equation}\label{t2}||\sigma(t,\omega,x)||^2+2\langle x,b(t,\omega,x)\rangle \leq
f(t,\omega)(\gamma(|x|^2)+1),\  |x|\ge K
 \end{equation}
with $f$ being a progressively measurable and non negative function satisfying
$$\mathbb{E} \int_0^t f(s,\cdot)ds<\infty,\ t\geq 0.$$
Then, the solution of equation $(\ref{sdelan})$ is global, namely, the lifetime
$$\zeta:=\inf\{t>0,|X_t|=\infty\}=+\infty.$$
\end{Theorem}

We now give an example to support our conditions (\ref{t3}) and (\ref{t2}). However, our example
does not fulfill the conditions ($H1$) and ($H2$) in \cite{Fang}, respectively, neither does for the conditions
of Theorem 4 in \cite{Watanabe}. This then indicates that our conditions are indeed weaker than those
known conditions.

\noindent\textbf{Example}\quad When the equation reduces to
time indendent case, let $m=d$, and
$\sigma(x)=diag(\sigma_1(x),\cdots,\sigma_d(x))$,
$b(x)=(b_1(x),\cdots,b_d(x))^T$, where
$\sigma_i(x)=x_i^{\frac{2}{3}},$ $b_i(x)=-x_i^{\frac{1}{3}}$, here
"$T$" denotes the transpose of the vector or a matrix. Then it is easy to see
that neither conditions ($H1$), ($H2$) in \cite{Fang} nor conditions in Theorem 4 in \cite{Watanabe}
hold since in this case the coefficients $\sigma$ and $b$ are
H$\ddot{\textrm{o}}$lder continuous with order $\frac{2}{3}$ and
$\frac{1}{3}$, respectively. But our conditions (\ref{t3}) and (\ref{t2}) are
satisfied for such defined $\sigma$ and $b$. Indeed,

\begin{eqnarray*}&& ||\sigma(x)-\sigma(y)||^2+2\langle x-y,b(x)-b(y)\rangle \\
&=&\sum_{i=1}^d(x_i^\frac{2}{3}-y_i^\frac{2}{3})^2-2\sum_{i=1}^d(x_i-y_i)(x_i^\frac{1}{3}-y_i^\frac{1}{3})\\&
=&\sum_{i=1}^d(x_i^\frac{1}{3}-y_i^\frac{1}{3})^2[(x_i^\frac{1}{3}+y_i^\frac{1}{3})^2-2(x_i^\frac{2}{3}+x_i^\frac{1}{3}y_i^\frac{1}{3}+y_i^\frac{2}{3})]\\&
=&-\sum_{i=1}^d(x_i^\frac{1}{3}-y_i^\frac{1}{3})^2(x_i^\frac{2}{3}+y_i^\frac{2}{3}) \\
&\le &\eta_R(|x-y|^2)\end{eqnarray*}
where we have used the fact that the second last line is clearly non-positive in the last derivation. Similarly, we have
$$||\sigma(x)||^2+2\langle x,b(x)\rangle
=\sum_{i=1}^d(x_i^\frac{2}{3})^2-2\sum_{i=1}^dx_ix_i^\frac{1}{3}
=-\sum_{i=1}^dx_i^\frac{4}{3}\le\gamma(|x|^2).
$$
This shows that our conditions (\ref{t3}) and (\ref{t2}) are fulfilled. Hence,
by Theorem \ref{dingli1} and Theorem \ref{dingli2}, there is a unique strong global solution of
the SDE (\ref{sdelan}).

We would like to point out here that in the above example the condition are given for  the two coefficients $b$ and
$\sigma$ jointly, which guarantee  the H\"older continuity condition. For instance, putting the drift coefficient $b$ to
be zero, the diffusion coefficient $\sigma$ in our example is clearly not H\"older continuous.

\begin{Remark}In $\cite{Lan}$, the first named author studied pathwise uniqueness and
non explosion properties of the solution of equation $(\ref{sdelan})$. We have
generalized Fang and Zhang's results to a more general case. But it is not clear at that time whether
there is a solution of equation $(\ref{sdelan})$ under the given condition. However, we can get the same
conclusion when $(\ref{t2})$ and $(\ref{t3})$ hold. Moreover, we can prove that there is really a
unique strong global solution of equation $(\ref{sdelan})$ if $(\ref{t2})$ and $(\ref{t3})$ hold.
\end{Remark}

In \cite{Hof},  Hofmanov\textrm{\'{a}} and Seidler proved the following result.
Assume $b$ and $\sigma$ are Borel functions such that $b(t,\cdot)$ and $\sigma(t,\cdot)$ are continuous for any $t\in [0,T]$ and the linear growth hypothesis is satisfied, that is

\[\exists K_*<\infty,\ \forall t\in [0,T],\ \forall x\in\mathbb{R}^m,\ ||b(t,x)||\vee||\sigma(t,x)||\le K_*(1+||x||).\]

Let $\nu$ be a Borel probability measure on $\mathbb{R}^m.$ Then there exists a weak solution to the problem

\[dX=b(t,X)dt+\sigma(t,X)dW, X(0)\sim\nu\]
where $W$ is a standard Brownian motion.

\begin{Remark}
Let us comment our conditions with those of Hofmanov$\acute{a}$ and Seidler in \cite{Hof}.
On one hand, our coercive condition \ref{t2} on the coefficients $b$ and $\sigma$ is weaker than the
corresponding linear growth hypothesis given by Hofmanov$\acute{a}$ and Seidler in \cite{Hof},
as we don't need the linear growth hypothesis here. In \cite{Hof}, however, the authors has proved the existence
of weak solutions under the spatial continuity of the coefficients plus linear growth hypothesis. In our paper,
we show the existence of a unique strong solution under the spatial continuity of the coefficients together
with a locally weak monotonicity condition.
\end{Remark}

Let us give some comparison of our conditions with those existing in the literature. There are many works dealing
with the existence and uniqueness of SDEs. Stroock and Varadhan (see \cite{Stroock}) proved the weak existence and uniqueness of equation (\ref{sdelan}) by using martingale problem method when $\sigma$ is bounded continuous and uniformly elliptic and $b$ is bounded and measurable, and both the coefficients are independent with $t$ and $\omega.$

In Watanabe and Yamada \cite{Watanabe,Yamada1},
the authors gave sufficient conditions on $\sigma$ and $b$ for the strong
uniqueness and existence of stochastic differential equation
$$dX_t=\sigma(t,X_t)dB_t+b(t,X_t).$$ If we take $\eta_R(x)=R(\rho^2(\sqrt{x})+\sqrt{x}\bar{\rho}(\sqrt{x}))$, where
$\rho$ and $\bar{\rho}$ are the same as that of \cite{Watanabe}, it's obvious that
condition (\ref{t3}) holds for such defined $\eta_R.$ Note that we
don't need the concave condition on $\rho$ and $\bar{\rho}$.

More recently, Fang and Zhang \cite{Fang} gave the sufficient conditions on
$\sigma$ and $b$ under which the degenerated time-homogeneous equation of (\ref{sdelan}) has no
explosion, pathwise uniqueness and non confluence. They proved a special case in \cite{Fang1} for non explosion and
pathwise uniqueness of the equation. Since $\sigma,\ b$ are both continuous, according to Ikeda and
Watanabe \cite{Ikeda} Chapter IV, Theorem 2.3, the solution does
exist under Fang and Zhang's conditions. By taking
$\eta_R(x)=Rxr(x), \gamma(x)=x\rho(x)+1$, where $r(x)$  and $\rho$ are the same as that of conditions
$(H1)$ and $(H2)$ in Fang and Zhang \cite{Fang}, our condition (\ref{t3}) and (\ref{t2}) are
satisfied. By Theorem \ref{dingli1} and Theorem \ref{dingli2}, there is a unique strong
solution of equation (\ref{sdelan}), which is non explosive. And the above example shows that
$(H1)$ and $(H2)$ in \cite{Fang} don't hold, but (\ref{t3}) and (\ref{t2}) hold. Moreover,
they must assume that the control function be differentiable to make
sure the Gronwall lemma can be used, but we can drop the
differentiability condition by using a new test function.
$\eta_R(x)=Rx\log(1/x),(x<1)$ is a typical example for our $\eta_R$.

In \cite{Rockner}, Pr$\acute{\textrm{e}}$v$\ddot{\textrm{o}}$t and
R$\ddot{\textrm{o}}$ckner (see also \cite{Krylov}) proved that when
$\sigma, b$ satisfy (\ref{t1}) and the so called weak coercivity and
local weak monotonicity, there exists a unique (up to
P-indistinguishability) solution to the stochastic differential
equation (\ref{sdelan}).  Both \cite{Krylov,Rockner} had to use the linearity of the control
function to prove that the approximation sequence $X^{(n)}(t)$ is
uniformly convergent. Since $\eta_R$ in our condition (\ref{t3}) may
not be linear function, and there is no weak coercivity, their
method can not be used in our case either. Krylov and R\"{o}ckner \cite{Krylov1} proved existence and
uniqueness of strong solutions to stochastic equations in domains
$G\subset \mathbb{R}^d$ with singular time dependent drift $b$ up to
an explosion time, but they must assume the unit diffusion.

Recently, Shao, Wang and Yuan \cite{Wang} proved that there exists a
unique non explosive solution for (\ref{sdelan}) when the
coefficients satisfy certain global conditions. Actually, their
assumptions guarantee the local coercivity condition in
\cite{Rockner,Krylov}. Thus, the condition is too strong in certain sense for the existence of (\ref{sdelan}).

We now turn to the moment estimation for the following Markovian type stochastic differential equation

\begin{equation}\label{sdelan2}dX_t=\sigma(t,X_t)dB_t+b(t,X_t)dt.\
X_0=x_0, \end{equation}

In \cite{HMS}, the authors get the upper bound of $p$th moment of maximum
process $\sup_{0\le s\le t}|X_t|$ for $X_t$ being the unique strong solution of the time-homogeneous SDE.
When the diffusion term is Lipschitz
continuous, and the drift term satisfies one-sided Lipschitz condition, they
prove that for $p\ge 2$, there exists $C(p,t)>0$ such that
$$\mathbb{E}(\sup_{0\le s\le t}|X_t|^p)\le C(p,t)(1+|x_0|^p).$$

Starting with that $X_t$ is the unique solution of our equation (\ref{sdelan2}),
we will investigate the $p$th moment of the solution under more general and weaker condition.

\begin{Theorem}\label{dingli4} Assume that the coefficients $\sigma$ and $b$ satisfy
\begin{equation}\label{t5}\Big(||\sigma(t,x)||^2+2\langle x,b(t,x)\rangle\Big)\vee |\sigma^T(t,x)x|^2
\leq f(t)(|x|^2+1).
\end{equation}
Let $p>2$ be fixed arbitrarily. We have the following

(i) If
$$0\le f\in L_{loc}^p(\mathbb{R}_+):=\{f:\mathbb{R}_+\rightarrow\mathbb{R}_+,\int_0^tf(s)^pds<\infty,\forall t>0\},$$
then we have
\[\mathbb{E}(\sup_{0\le s\le t}|X_s|^p)\le  A\exp\{B\int_0^tf(s)^\frac{p}{2}ds+C\int_0^tf^p(s)ds\}\]
where $A, B, C$ are only dependent of $p,$ $t$ and function $f$.

(ii) If $0\le f\in L_{loc}^{\frac{2p}{p-2}}(\mathbb{R}_+)$, then for any fixed $t>0,$
\[\mathbb{E}(\sup_{0\le s\le t}|X_s|^p)\le A_1e^{B_1t}\]
where $A,B$ are still only dependent of $p,$ $t$ and function $f$.
\end{Theorem}

\begin{Remark}
We would like to point out that when $p$ is close to 2, $L_{loc}^\frac{2p}{p-2}$ is
close to $L_{loc}^\infty,$ so $f\in L_{loc}^\frac{2p}{p-2}$ is essentially that $f$
is bounded. But $f\in L_{loc}^p$ is not necessarily bounded. On the other hand, for
$p$ sufficiently large, $f\in L_{loc}^\frac{2p}{p-2}$ is near that $f\in L_{loc}^2$,
but $f\in L_{loc}^p$ is near essentially bounded.
\end{Remark}

Finally, let us consider the property of non confluence of the following time-homogeneous SDE

\begin{equation}\label{sdelan1}dX_t=\sigma(X_t)dB_t+b(X_t)dt,\
X_0=x_0, \end{equation}

We say that the solution $X_t$ of equation (\ref{sdelan1}) has
non confluence, if for all $x_0\neq y_0,$
$$P(X_t(x_0)\neq X_t(y_0),\ \forall t>0)=1. $$

Such kind of non confluence property was studied by Emery in an
early work \cite{Emery} for general stochastic differential
equations under Lipschitzian conditions, and by Yamada and Ogura for
non-Lipschitz case in \cite{Yamada}. However the mixing condition
imposed in \cite{Yamada} for coefficients $\sigma$ and $b$ is
difficult to be checked and not natural.

Fix $R>0$ arbitrarily, let $\gamma_R:[0,1)\rightarrow\mathbb{R}_+$ be a
differentiable function on a very small interval
$[0,\varepsilon_0)(0<\varepsilon_0\le c_0)$  such that
$$\gamma_R(0)=0,\ \int_{0+}\frac{dx}{\gamma_R(x)}=\infty$$
and
$$\frac{x(\gamma'_R(x)+1)}{\gamma_R(x)}\le K,\ \forall |x|\le c_0$$
for some constant  $K> \frac{1}{2}$ which is independent of $x$ and $R$.

Then by using a new test function, we can show the following result

\begin{Theorem}\label{dingli3}  If the coefficients $\sigma$ and $b$ satisfy that for
$|x|\vee|y|\le R,$
\begin{equation}\label{t4}||\sigma(x)-\sigma(y)||^2-\frac{2}{2K-1}\langle x-y,b(x)-b(y)\rangle
\leq \gamma_R(|x-y|^2),\ |x-y|\le c_0,
\end{equation}
then the solution $X_t$ of equation (\ref{sdelan1}) has non confluence property before life time $\zeta$. That is,
for any $x_0\neq y_0,$ the conditional probability
$$P(\{\omega\in\Omega:X_t(x_0,\omega)\neq X_t(y_0,\omega),t>0\}|\{\omega\in\Omega: t<\zeta(\omega)\})=1.$$\end{Theorem}

For the case of one dimension, the non confluence corresponds to stochastic monotonicity, that is, if $x_0\le y_0$,
 then $X_t(x_0)\le X_t(y_0)$ for all $t\ge 0$ almost surely. Since condition (\ref{t4}) naturally holds when $\sigma$ and $b$ are locally Lipschitzian continuous, then we have the following

\begin{Corollary}
Suppose d=m=1. If $\sigma$ and $b$ are locally Lipschitzian continuous, then the solution $X_t$ of equation (\ref{sdelan1}) is stochastic monotone before life time.
\end{Corollary}

The rest of the paper is organized as follows. We show there exists uniformly convergent subsequence of Euler approximation which is our Lemma \ref{l1} in Section 2. Then we prove that the limit process is a solution of the equation (\ref{sdelan}) in Section 3. In Section 4, we will prove non explosion result which is our Theorem \ref{dingli2}. Then in
Section 5 we will get the upper bound of the $p$th moment of the maximum process. Finally, we prove non confluence
of the solution of time-homogeneous SDE in Section 6.

\section{Uniform convergence of Euler approximaiton}\vskip0.2in

\noindent \textbf{Proof of Lemma \ref{l1}}  If we denote $\kappa(n,t):=[tn]/n$, equation (\ref{t8}) is
equivalent to

\begin{equation}\label{bj} X^{(n)}(t)=x_0+\int_0^tb(s,\cdot,X^{(n)}(s)
+p^{(n)}(s))ds+\int_0^t\sigma(s,\cdot,X^{(n)}(s)+p^{(n)}(s))dB_s\end{equation}
where
\begin{equation}\label{pnt}
\aligned p^{(n)}(t)&:=X^{(n)}(\kappa(n,t))-X^{(n)}(t)\\&
=-\int_{\kappa(n,t)}^tb(s,\cdot,X^{(n)}(\kappa(n,s))ds\\&
-\int_{\kappa(n,t)}^t\sigma(s,\cdot,X^{(n)}(\kappa(n,s))dB_s,t\in[0,\infty).\endaligned\end{equation}

In what follows, we want to prove that there is a subsequence of $X^{(n)}(t) $ converges to some
process $X(t)$. Define $$ \tau^{(n)}(R):=\inf\{t>0,|X^{(n)}(t)|\geq R\},$$ then by
 the definition of $p^{(n)}(t)$ in (\ref{pnt}),
$$\aligned -\langle e_i,p^{(n)}(t)\rangle= &\int_{\kappa(n,t)}^t\langle e_i,b(s,\cdot,X^{(n)}(\kappa(n,s))\rangle ds
+\int_{\kappa(n,t)}^t\langle
e_i,\sigma(s,\cdot,X^{(n)}(\kappa(n,s))\rangle dB_s\endaligned$$
where $ e_i,1\le i\le d$, is the canonical basis of $ \mathbb{R}^d.$
It follows that

\begin{equation}\label{pnt1}\aligned &P(|\langle e_i,p^{(n)}(t)\rangle|\ge 2\varepsilon,t\le\tau^{(n)}(R))\\ &\le P(\int_{\kappa(n,t)}^t\sup_{|x|\le R}|b(s,x)|ds\ge\varepsilon)
\\&\quad+P(\sup_{\tilde{t}\in[0,t]}|\int_0^{\tilde{t}\wedge \tau^{(n)}(R)}1_{[\kappa(n,t),\infty)}(s)\langle e_i,\sigma(s,\cdot,X^{(n)}(\kappa(n,s))\rangle dB_s|\ge\varepsilon)\\&
\le P(\int_{\kappa(n,t)}^t\sup_{|x|\le
R}|b(s,x)|ds\ge\varepsilon)+\frac{1}{\varepsilon^2}\mathbb{E}\int_{\kappa(n,t)}^{t\wedge\tau^{(n)}(R)}\sup_{|x|\le
R}||\sigma(s,x)||^2ds.\endaligned\end{equation}
The last inequality is derived  by utilising martingale inequality. Then for any
fixed $t\in [0,\infty)$ and $R$, we have

\begin{equation}\label{pntjx} 1_{\{t\le
\tau^{(n)}(R)\}}p^{(n)}(t)\stackrel{P}{\longrightarrow } 0,\
n\rightarrow \infty. \end{equation}
Hence, there exists a subsequence of $\{n\}$ depending on $R$ and $t$ (which is also denoted as $\{n\}$ for the
sake of simplicity) such that

\begin{equation}\label{pntas} 1_{\{t\le
\tau^{(n)}(R)\}}p^{(n)}(t)\stackrel{a.s.}{\longrightarrow }
0,n\rightarrow \infty. \end{equation}

Now let $Y_t=X^{(n)}(t)-X^{(m)}(t)$, $\xi_t=|Y_t|^2 .$ Define the
following test function:

\begin{equation}\label{pntas1}\varphi_\delta(x)=\int_0^x\frac{ds}{\eta_R(s)+\delta}. \end{equation}

It's obvious that for any $ 0<x<\varepsilon_0, $ when
$\delta\downarrow 0 $,

\begin{equation}\label{phi}
\varphi_\delta(x)\uparrow\varphi_0(x)=\int_0^x\frac{ds}{\eta_R(s)}=\infty.\end{equation}

Since

\begin{equation}\label{phijx}\frac{\varphi_\delta(x_2)-\varphi_\delta(x_1)}{x_2-x_1}\ge
\frac{\varphi_\delta(x_3)-\varphi_\delta(x_2)}{x_3-x_2},\quad
0<x_1<x_2<x_3<\varepsilon_0,
\end{equation}
$\varphi_\delta(x)$ is a concave function on the interval $
[0,\varepsilon_0)$. Note here that since $\lim_{x\rightarrow 0}\varphi'_\delta(x)=\frac{1}{\delta}$,
there is a concave extension of $\varphi_\delta(x)$ on the real line.
Let $ \bar{\varphi}_\delta $ be a concave
function on $\mathbb{R}$, and satisfy $
\bar{\varphi}_\delta(x)=\varphi_\delta(x),x\in [0,\varepsilon_0). $
Then the second order derivative of $ \bar{\varphi}_\delta $ in the
sense of distributions $ \bar{\varphi}{''}_\delta $ is a non
positive Radon measure(see \cite{Revuz} Appendix 3). Let

$$
\tau_{n,m}:=\inf\{t>0,\xi_t\geq\varepsilon_0\}.$$

By using $ \textrm{It}\hat{\textrm{o}}$-$\textrm{Tanaka's}$ formula,
we have
\begin{eqnarray*}
&& \bar{\varphi}_{\delta}(\xi_{t\wedge\tau_{n,m}})\\
&=&2\int_0^{t\wedge\tau_{n,m}}\bar{\varphi'_{\delta}}(\xi_s)
\langle Y_s,(\sigma(s,\omega,X^{(n)}(s)+p^{(n)}(s))-\sigma(s,\omega,X^{(m)}(s)+p^{(m)}(s))dB_s\rangle\\
&&\quad+2\int_0^{t\wedge\tau_{n,m}}\bar{\varphi}'_{\delta}(\xi_s)
\langle
Y_s,b(s,\omega,X^{(n)}(s)+p^{(n)}(s))-b(s,\omega,X^{(m)}(s)+p^{(m)}(s))\rangle
ds\\
&&\quad+\int_0^{t\wedge\tau_{n,m}}\bar{\varphi}'_{\delta}(\xi_s)
|| \sigma(s,\omega,X^{(n)}(s)+p^{(n)}(s)))-\sigma(s,\omega,X^{(m)}(s)+p^{(m)}(s)))||^2ds\\
&&\quad
+\frac{1}{2}\int_{\mathbb{R}}L_{t\wedge\tau_{n,m}}^a(\xi)\bar{\varphi}{''}_{\delta}(da)
\\&\le& 2\int_0^{t\wedge\tau_{n,m}}\varphi'_{\delta}(\xi_s)
\langle Y_s,(\sigma(s,\omega,X^{(n)}(s)+p^{(n)}(s)))-\sigma(s,\omega,X^{(m)}(s)+p^{(m)}(s))dB_s\rangle\\
&&\quad+2\int_0^{t\wedge\tau_{n,m}}\varphi'_{\delta}(\xi_s) \langle
Y_s,b(s,\omega,X^{(n)}(s)+p^{(n)}(s))-b(s,\omega,X^{(m)}(s)+p^{(m)}(s)))\rangle
ds\\
&&\quad+\int_0^{t\wedge\tau_{n,m}}\varphi'_{\delta}(\xi_s) ||
\sigma(s,\omega,X^{(n)}(s)+p^{(n)}(s)))-\sigma(s,\omega,X^{(m)}(s)+p^{(m)}(s))||^2ds.
\end{eqnarray*}

The above inequality holds because the second derivative $
\bar{\varphi}{''}_\delta $ in the sense of distributions is a non
positive Radon measure and the local time is always non positive.
Since
$\varphi_{\delta}(\xi_{t\wedge\tau})=\bar{\varphi}_{\delta}(\xi_{t\wedge\tau}),$
taking expectation on both sides, we have
\begin{equation}\label{qiwang}
\aligned \mathbb{E}\ \varphi_{\delta}(\xi_{t\wedge\tau_{n,m}})&
\!\!\leq\!\! 2\
\mathbb{E}\!\!\int_0^{t\wedge\tau_{n,m}}\!\!\varphi'_{\delta}(\xi_s)
\langle
Y_s,b(s,\omega,X^{(n)}(s)\!\!+\!\!p^{(n)}(s))\!\!-\!\!b(s,\omega,X^{(m)}(s)\!\!+\!\!p^{(m)}(s))\rangle
ds\\
&\quad+\mathbb{E}\!\!\int_0^{t\wedge\tau_{n,m}}\!\!\!\!\varphi'_{\delta}(\xi_s)
||
\sigma(s,\omega,X^{(n)}(s)\!\!+\!\!p^{(n)}(s)))\!\!-\!\!\sigma(s,\omega,X^{(m)}(s)\!\!+\!\!p^{(m)}(s))||^2ds.
\endaligned\end{equation}
By the definition of $p^{(n)}(t) $ given in (\ref{pnt}), we know that the above
quantity could be dominated by
\begin{equation}\label{qiwang1}
\aligned
&2\mathbb{E}\!\!\int_0^{t\wedge\tau_{n,m}}\!\!\varphi'_{\delta}(\xi_s)
\langle
p^{(m)}(s)-p^{(n)}(s),b(s,\omega,X^{(n)}(s)+p^{(n)}(s))-b(s,\omega,X^{(m)}(s)+p^{(m)}(s))\rangle
ds\\
&\quad+\mathbb{E}\!\!\int_0^{t\wedge\tau_{n,m}}\!\!\varphi'_{\delta}(\xi_s)g(s,\omega)
\eta_R(|X^{(n)}(s)+p^{(n)}(s)-X^{(m)}(s)-p^{(m)}(s)|^2)ds.
\endaligned\end{equation}

Denote \begin{equation}\label{hs}H(s):=|\varphi'_{\delta}(\xi_s)
\langle
p^{(m)}(s)\!-\!p^{(n)}(s),b(s,\omega,X^{(n)}(s)\!+\!p^{(n)}(s))\!-\!b(s,\omega,X^{(m)}(s)\!+\!p^{(m)}(s))\rangle|,\end{equation}
and
\begin{equation}\label{gs}G(s):=\varphi'_{\delta}(\xi_s)\eta_R(|X^{(n)}(s)+p^{(n)}(s)-X^{(m)}(s)-p^{(m)}(s)|^2).\end{equation}

Then
$$
H(s)\le 2\sup_{|x|\le
R}|b(s,\omega,x)|\times\frac{1}{\eta_R(\xi_s)+\delta}|p^{(n)}(s)-p^{(m)}(s)|,\quad
s\le \tau^{(n,m)}(R)
$$
where $$ \tau^{(n,m)}(R):=\tau^{(n)}(R)\wedge\tau^{(m)}(R).$$

Denote

\begin{equation}T_{n,m}(R):=\tau_{n,m}\wedge\tau^{(n,m)}(R).\end{equation}

Then, by (\ref{pntas}), for fixed $\delta,$ let $m,n$ be large
enough(dependent on $\delta$), we have
$$\aligned
\frac{1}{\eta_R(\xi_s)+\delta}|p^{(n)}(s)-p^{(m)}(s)|&\le\frac{1}{\delta}|p^{(n)}(s)-p^{(m)}(s)|\\&\le
1,\quad s\in[0,t\wedge T_{n,m}(R)), a.s..
\endaligned
$$
So
\begin{equation}\label{hs1}H(s)\le 2\sup_{|x|\le R}|b(s,\omega,x)|,\quad
s\in[0,t\wedge T_{n,m}(R)), a.s..\end{equation}

Similarly, by the continuity of function $\eta_R$, for fixed $\delta$,
and $m,n$ be large enough(dependent on $\delta$), we have

\begin{equation}
\label{hs1a}\aligned
G(s)&=\frac{\eta_R(|X^{(n)}(s)-X^{(m)}(s)+p^{(n)}(s)-p^{(m)}(s)|^2)}{\eta_R(|X^{(n)}(s)-X^{(m)}(s)|^2)+\delta}\\&
\le 1+\frac{|\eta_R(|X^{(n)}(s)-X^{(m)}(s)+p^{(n)}(s)-p^{(m)}(s)|^2)
-\eta_R(|X^{(n)}(s)-X^{(m)}(s)|^2)|}{\eta_R(|X^{(n)}(s)-X^{(m)}(s)|^2)+\delta}\\&\le
2,\quad s\in[0,t\wedge T_{n,m}(R)),
a.s..\endaligned\end{equation}

Then by Fatou's lemma, we have

\begin{equation}\label{qiwang2}
\aligned
&\limsup_{m,n\rightarrow\infty}\mathbb{E}\!\!\int_0^{t\wedge\tau_{n,m}}\!\!\varphi'_{\delta}(\xi_s)g(s,\omega)
\eta_R(|X^{(n)}(s)+p^{(n)}(s)-X^{(m)}(s)-p^{(m)}(s)|^2)ds\\&
\le\mathbb{E}\!\!\int_0^{t}\!\!\limsup_{m,n\rightarrow\infty}\textbf{1}_{[0,\tau_{n,m}]}(s)\varphi'_{\delta}(\xi_s)g(s,\omega)
\eta_R(|X^{(n)}(s)+p^{(n)}(s)-X^{(m)}(s)-p^{(m)}(s)|^2)ds\\&\le
2\mathbb{E}\int_0^tg(s,\omega)ds.
\endaligned\end{equation}

In the last inequality above we used the fact that
$$\limsup_{m,n\rightarrow\infty}\varphi'_{\delta}(\xi_s)
\eta_R(|X^{(n)}(s)+p^{(n)}(s)-X^{(m)}(s)-p^{(m)}(s)|^2)ds$$
is bounded almost surely.
Hence
$$\mathbb{E}\int_0^{t\wedge\tau_{n,m}}\varphi'_{\delta}(\xi_s)g(s,\omega)
\eta_R(|X^{(n)}(s)+p^{(n)}(s)-X^{(m)}(s)-p^{(m)}(s)|^2)ds$$ is
uniformly bounded when $m,n$ are large. Now obviously the choice of
$m$ and $n$ is independent of sample points. An argument with
similar spirit may also work for the first term of right hand side
in (\ref{qiwang1}). Hence

\begin{equation}
\label{qwjx}\limsup_{n,m\rightarrow\infty}\
\mathbb{E}\varphi_{\delta}(\xi_{t\wedge T_{n,m}(R)})\le
C(t)<\infty\end{equation}
where $C(t)$ is a positive constant
independent of $\delta$ and $R$, $m(\delta)$ and $n(\delta)$.

Next, suppose that there are subsequences $\{m_k\}, \{n_k\}$ such that

$$
\lim_{k\rightarrow\infty}\xi_{t\wedge T_{(n_k,m_k)}(R)}=\limsup_{n,m\rightarrow\infty}\
\xi_{t\wedge T_{n,m}(R)}=\xi_0,\quad
\textrm{as}\ k\rightarrow\infty.$$

Since

$$|\varphi_{\delta}(\xi_{t\wedge T_{n,m}(R)})|\le\frac{\varepsilon_0}{\delta}$$

for any $n,m$ when $\delta$ is fixed, by using the dominated
convergence theorem, we have

$$ \aligned\mathbb{E}\varphi_{\delta}(\limsup_{n,m\rightarrow\infty}\
\xi_{t\wedge T_{n,m}(R)})
&=\lim_{k\rightarrow\infty}\mathbb{E}\varphi_{\delta}(\xi_{t\wedge T_{(n_k,m_k)}(R)})
\\&\le\limsup_{n,m\rightarrow\infty}\
\mathbb{E}\varphi_{\delta}(\xi_{t\wedge T_{n,m}(R)})\le
C(t)<\infty.\endaligned
$$

Then let $\delta\downarrow 0 $. Since $C(t)$ is independent of
$\delta,$ we have
$$
  \mathbb{E}\int_0^{\limsup\limits_{n,m\rightarrow\infty} \xi_{t\wedge T_{n,m}(R)}} \frac{1}{\eta_R(s)}ds\le C(t)<\infty.
$$

By using $\int_{0+}\frac{ds}{\eta_R(s)}=\infty,$ it follows that for any fixed $t>0,$

\begin{equation}\label{gljx1}
P(\limsup_{n,m\rightarrow\infty}\ \xi_{t\wedge T_{n,m}(R)}=0)=1.\end{equation}

Now by Fatou's lemma we have
\begin{equation} \label{qwjx1}
 0\le\limsup_{n,m\rightarrow\infty}
 \ \mathbb{E}\ \xi_{t\wedge T_{n,m}(R)}\le\mathbb{E}\ \limsup_{n,m\rightarrow\infty}
 \xi_{t\wedge T_{n,m}(R)}=0,
\end{equation}
that is, $\lim\limits_{n,m\rightarrow\infty} \ \mathbb{E}\ \xi_{t\wedge T_{n,m}(R)}$ does exist and the limit is 0. Therefore,
\begin{equation}
\label{budengshi1}\aligned P(\sup_{s\le
t\wedge T_{n,m}(R)}|X^{(n)}(s)-X^{(m)}(s)|^2\ge\varepsilon)\le\frac{\mathbb{E}\
\xi_{t\wedge T_{n,m}(R)}}{\varepsilon}\rightarrow0.\endaligned\end{equation}

 We can now select a subsequence, which will again be denoted by
 $X^{(n)}$ such that

\begin{equation}
\label{budengshi2}
P(\sup_{s\let\wedge T_{n,m}(R)}|X^{(n)}(s)-X^{(m)}(s)|^2\ge2^{-m\wedge n}) \le 2^{-m\wedge n}.
\end{equation}
Since $t$ is arbitrary, we have

\begin{equation}
\label{budengshi3}\aligned P(\sup_{s\le
T_{n,m}(R)}|X^{(n)}(s)-X^{(m)}(s)|^2\ge
2^{-m\wedge n}) \le 2^{-m\wedge n}.\endaligned\end{equation}

Denote
$$\tau_R:=\liminf_{n,m\to\infty}T_{n,m}(R).$$

Due to (\ref{budengshi3}) there is a subsequence $X^{(n_k)}$ which is convergent to an $\mathscr{F}_t$ adapted process
$X$ defined in $[0,\tau_R]$ $\mathbb{P}$-almost surely in $C([0,\tau_R];\mathbb{R}^d)$.

Let $R\rightarrow\infty$, we obtain that there is a subsequence of $X^{(n_{k})}$ (still denotes $X^{(n_{k})}$ for simplicity)
such that for any fixed $T>0,$
\begin{equation}
\label{shoulian1}\sup_{t\le\tau(T)}|X^{(n_k)}(t)-X(t)|
 \stackrel{a.s.}{\longrightarrow}0,\ k\rightarrow\infty \end{equation}
where $$\tau(T):=\liminf_{R\rightarrow\infty} \tau_R\wedge T.$$
The proof is complete. $\square$

\section{The limit process is a solution of equation (\ref{sdelan})}

\noindent \textbf{Proof of Lemma \ref{l2}} We are now going to prove that the limit process $X(t)$ is a
solution of SDE (\ref{sdelan}). By
(\ref{shoulian1}) we only need to prove that there is a subsequence of the right hand side of
(\ref{bj}) converges to

$$X_0+\int_0^tb(s,\cdot,X(s))ds+\int_0^t\sigma(s,\cdot,X(s))dB_s,\quad t\le\tau(T).$$

Since the convergence in (\ref{shoulian1}) is uniform, by
equicontinuity we have

\begin{equation}
\label{tau}\sup_{t\le
\tau(T)}|X^{(n_k)}(\kappa(n_k,t))-X(t)|
\stackrel{a.s.}{\longrightarrow}0,\ k\rightarrow\infty.
\end{equation}

In particular, for $S(t):=\sup_
{k\in\mathbb{N}}|X^{(n_k)}(\kappa(n_k,t))|,$

\begin{equation}
\label{tau1}\sup_{t\le\tau(T)}S(t)<\infty,\
P-a.s..\end{equation}

For $R\in[0,\infty),$ define the $(\mathscr{F}_t)$-stopping time

$$\tilde\tau(R,T):=\inf\{t\ge 0,S(t)>R\}\wedge\tau(T).$$

By the continuity of $b$ in $x\in\mathbb{R}^d$ and by local
integrability condition (\ref{t1})

\begin{equation}
\label{jifen}\lim_{k\rightarrow\infty}\int_0^tb(s,\cdot,X^{(n_k)}(\kappa(n_k,s))ds
=\int_0^tb(s,\cdot,X(s))ds,\quad P-a.s.\ \textrm{on}\
\{t\le\tilde\tau(R,T)\}.\end{equation}

For the stochastic integrals part we construct another sequence of
stopping times. For $R,N\ge 0,$ define

$$\tau_N(R,T):=\inf\{t\ge 0,\int_0^T\sup_{|x|\le
R}||\sigma(s,\cdot,x)||^2ds\ge N\}\wedge\tilde\tau(R,T).$$

Now by the continuity of $\sigma$ in $x\in\mathbb{R}^d$, (\ref{t1})
and Lebesgue's dominated convergence theorem

$$\lim_{k\rightarrow\infty}\mathbb{E}\int_0^{\tau_N(R,T)}||\sigma(s,\cdot,X^{(n_k)}(\kappa(n_k,s)))-\sigma(s,\cdot,X(s)||^2ds=0,$$
hence

\begin{equation}
\label{jifen1}P\textrm{-}\lim_{k\rightarrow\infty}\int_0^t\sigma(s,\cdot,X^{(n_k)}(\kappa(n_k,s)))dB_s
=\int_0^t\sigma(s,\cdot,X(s)dB_s\end{equation}
 on $t\le\tau_N(R,T).$

Thus, there exists a subsequence of $\{n_k\}$ (which is also denoted as $\{n_k\}$ for the sake of simplicity) such that

\begin{equation}
\label{jifen2}\lim_{k\rightarrow\infty}\int_0^t\sigma(s,\cdot,X^{(n_k)}(\kappa(n_k,s)))dB_s
=\int_0^t\sigma(s,\cdot,X(s)dB_s,\,  a.s.\end{equation}

\noindent on
$t\le\tau_N(R,T).$
By the integrability condition (\ref{t1}) for every
$\omega\in\Omega$ there exists $N(\omega)\ge 0$ such that
$\tau_N(R,T)=\tilde\tau(R,T)$ for all $N\ge N(\omega),$ so

$$\bigcup_{N\in \mathbb{N}}\{t\le\tau_N(R,T)\}=\{t\le\tilde\tau(R,T)\}.$$

Therefore, (\ref{jifen2}) holds on $t\le\tilde\tau(R,T)$ almost surely. But by the
definition of $\tilde\tau(R,T)$ we have

$$\lim_{R\rightarrow\infty}\tilde\tau(R,T)=\tau(T).$$
It follows that
\begin{equation}
\label{jifen3}\aligned
X(t)=X_0+\int_0^tb(s,\cdot,X(s))ds+\int_0^t\sigma(s,\cdot,X(s))dB_s,
\quad
t\le\tau(T).\endaligned
\end{equation}
So the SDE (\ref{sdelan}) has a solution $X(t)$ at least before $t\le\tau(T)$ for any $T>0$. $\square$

\noindent \textbf{Proof of Theorem \ref{dingli1}}: By the same method that we have already used in \cite{Lan} and
cut by a stopping time $\tau_R=\inf\{t>0,|X_t|\vee|Y_t|\ge R\}$, where $X_t, Y_t$ are two solutions 
of SDE (\ref{sdelan}) with same initial value $x_0$. We can prove that SDE (\ref{sdelan})has pathwise 
uniqueness when the coefficients $\sigma,\ b $ satisfy condition (\ref{t3}). So by Lemma \ref{l1} 
and Lemma \ref{l2} there exists a unique strong solution of equation (\ref{sdelan}) with life time. $\square$

\section{Non explosion of the solution}

\textbf{Proof of Theorem \ref{dingli2}} First we show that
$$\liminf_{R\rightarrow\infty}\tau_R=\zeta, \, \textrm{a.s..}$$
By the definition of $\tau^{(n,m)}(R)$, it follows that $\tau^{(n,m)}(R)\rightarrow\zeta$ as $m, n, R\rightarrow\infty$ subsequently,
according to (\ref{gljx1}), we have

\begin{equation}\label{xijixian}
\limsup_{k\rightarrow\infty}\ \xi_{t\wedge\tau_{n_k,m_k}\wedge\zeta}=0,\textrm{P-a.s.}.
\end{equation}

If
$P(\liminf_{k\rightarrow\infty}\tau_{n_k,m_k}<\zeta)>0,$
then there exists a subsequence $\{k_l\}$ such that

$$P(\lim\limits_{l\rightarrow\infty}\tau_{n_{k_l},m_{k_l}}<\zeta)>0.$$

Therefore, for $0<T<\zeta$ close to $\zeta$ enough, there exists sufficiently large
$l>0$ such that

$$P(\tau_{n_{k_l},m_{k_l}}\le T<\zeta)>0.$$

It follows that on $\{\tau_{n_{k_l},m_{k_l}}\le
T<\zeta\}$,

$$\xi_{T\wedge\tau_{n_{k_l},m_{k_l}}\wedge\zeta}=\xi_{\tau_{n_{k_l},m_{k_l}}}.$$

According to (\ref{xijixian}), $\xi_{T\wedge\tau_{n_{k_l},m_{k_l}}\wedge\zeta}$ is smaller than any
fixed positive $\varepsilon$ for sufficiently large $l,$ we have

$$\xi_{\tau_{n_{k_l},m_{k_l}}}< \varepsilon_0\ \textrm{as}\ l\rightarrow\infty$$
holds with positive probability. Which is absurd since $\xi_{\tau_{n,m}}\equiv \varepsilon_0>0$ by
definition. So $\liminf_{k\rightarrow\infty}\tau_{n_k,m_k}\ge\zeta, \textrm{a.s.}$. By definition of $\tau_R,$ we have

$$\liminf_{R\rightarrow\infty}\tau_R=\liminf_{R\rightarrow\infty}\liminf_{n,m\rightarrow\infty}\tau_{n,m}\wedge\tau^{(n,m)}(R)=\zeta$$

Now we only need to show $\zeta=\infty$.

Define

$$\varphi(x):=\int_0^{x}\frac{ds}{\gamma(s)+1}.$$ Then
$\varphi$ is a concave function on $[0,\infty)$ since $\gamma$ is an increasing function and
$$\varphi'(x)=\frac{1}{\gamma(x)+1}.$$ As in the proof of Lemma \ref{l1}, we can extend $\varphi$ to a
function $\bar{\varphi}$ on $\mathbb{R}$ which is still concave.

Let $\tilde{\xi}_t:=|X_t|^2.$ Define

$$\hat{\tau}_R:=\inf\{t>0,\ |X_t|\geq R\},\quad R>0,$$
then $\hat{\tau}_R\uparrow\zeta$ as $R\uparrow\infty.$ By It\^{o}-Tanaka's formula, we
have

$$\aligned\varphi(\tilde{\xi}_{t\wedge\hat{\tau}_R})=\bar{\varphi}(\tilde{\xi}_{t\wedge\hat{\tau}_R})&
=\bar{\varphi}(\tilde{\xi}_0)+2\int_0^{t\wedge\hat{\tau}_R}\bar{\varphi}'(\tilde{\xi}_s)
\langle X_s,\sigma(s,\omega,X_s)dB_s\rangle\\
&\quad+2\int_0^{t\wedge\hat{\tau}_R}\bar{\varphi}'(\tilde{\xi}_s) \langle
X_s,b(s,\omega,X_s)\rangle ds
\\&\quad+\int_0^{t\wedge\hat{\tau}_R}\bar{\varphi}'(\tilde{\xi}_s)
||\sigma(s,\omega,X_s)||^2ds+\frac{1}{2}\int_R L^a_{t\wedge\hat{\tau}_R}(\tilde{\xi})\bar{\varphi}{''}(da)
\\&\le\varphi(\tilde{\xi}_0)+2\int_0^{t\wedge\hat{\tau}_R}\varphi'(\tilde{\xi}_s)
\langle X_s,\sigma(s,\omega,X_s)dB_s\rangle\\
&\quad+\int_0^{t\wedge\hat{\tau}_R}\varphi'(\tilde{\xi}_s)\big(2\langle
X_s,b(s,\omega,X_s)\rangle +
||\sigma(s,\omega,X_s)||^2\big)ds,
\endaligned $$
The inequality holds because $\bar{\varphi}{''}$ in the
sense of distribution is a non positive Radon measure and the local time
is non positive (see \cite{Revuz}, Appendix). Furthermore, taking expectation on both sides, we have

$$\aligned
\mathbb{E}\varphi(\tilde{\xi}_{t\wedge\hat{\tau}_R})&\le\varphi(\tilde{\xi}_0)+
\mathbb{E}\big(\int_0^{t\wedge\hat{\tau}_R}\varphi'(\tilde{\xi}_s) (2\langle
X_s,b(s,\omega,X_s)\rangle+
||\sigma(s,\omega,X_s)||^2)ds\big)\\&
\le\varphi(\tilde{\xi}_0)+\mathbb{E}\int_0^tf(s,\omega)ds
=:C_t<\infty.\endaligned$$
Note that $C_t$ is independent of $R.$ Letting $R\uparrow\infty$, and using Fatou's lemma, we
get
$$\mathbb{E}\int_0^{\tilde{\xi}_{t\wedge\zeta}}\frac{dx}{\gamma(x)+1}\leq
C_t<\infty.$$ Now if $P(\zeta<\infty)>0,$ then for some
$T>0$, $P(\zeta<T)>0.$ Taking $t=T$ in (2.4), we get
$$P(\zeta<T)\varphi(\tilde{\xi}_{\zeta})=\mathbb{E}\big({\bf 1}_{\{\zeta< T\}}\varphi(\tilde{\xi}_{\zeta})\big)\leq
C_T,$$ which is impossible since $\varphi(\tilde{\xi}_{\zeta})=\infty.$ Thus, for any $T>0,$

\begin{equation}
\label{jifen3a}\aligned
X(t)=X_0+\int_0^tb(s,\cdot,X(s))ds+\int_0^t\sigma(s,\cdot,X(s))dB_s,
\quad
t\le T.\endaligned
\end{equation}

It follows that the solution is non-explosive (That is,
the lifetime $\zeta=\infty$). We complete the proof. $\square$

\section{Moment inequality of the solution}

\noindent\textbf{Proof of Theorem \ref{dingli4}}\quad Denote $Y_t:=\sup_{0\le s\le t}|X_s|$. By It\^{o}'s formula, we have

\begin{equation}|X_t|^2=x_0^2+\int_0^t(2\langle X_s,b(s,X_s)+||\sigma(s,X_s)||^2)ds+M_t\end{equation}
where $M_t:=2\int_0^t\langle X_s,\sigma(s,X_s)dB_s\rangle.$ Then

\begin{equation}Y_t^2\le |x_0|^2+\int_0^tf(s)(|Y_s|^2+1)ds+\sup_{0\le s\le t}|M_s|.\end{equation}
So there exists $C_p>0$ such that

\begin{equation}Y_t^p\le C_p\Big(|x_0|^p+\big(\int_0^tf(s)(|Y_s|^2+1)ds\big)^\frac{p}{2}+\sup_{0\le s\le t}|M_s|^\frac{p}{2}\Big).\end{equation}
Thus,

\begin{equation}\label{moment}\mathbb{E}(Y_t^p)\le C_p\Big(|x_0|^p+\mathbb{E}\big(\int_0^tf(s)(|Y_s|^2+1)ds\big)^\frac{p}{2}+\mathbb{E}(\sup_{0\le s\le t}|M_s|^\frac{p}{2})\Big).\end{equation}
By Burkholder-Davis-Gundy inequality, there exists $C'_p>0$ such that

\[\aligned
\mathbb{E}(\sup_{0\le s\le t}|M_s|^\frac{p}{2})&\le C'_p\mathbb{E}[(\int_0^t|\sigma^T(s,X_s)X_s|^2ds)^\frac{p}{4}]\\&
\le C'_p\mathbb{E}\Big((Y^2_t+1)^\frac{p}{4}\big(\int_0^tf(s)|\sigma^T(s,X_s)X_s|ds\big)^\frac{p}{4}\Big)\\&
\le C'_p\Big(\frac{1}{2K}\mathbb{E}((Y^2_t+1)^\frac{p}{2})+\frac{K}{2}\mathbb{E}\big[\big(\int_0^tf^2(s)(|Y_s|^2+1)ds\big)^\frac{p}{2}\big]\Big).
\endaligned\]
Then

\[\aligned\mathbb{E}(Y_t^p)&\le C_p|x_0|^p+C_p\mathbb{E}\big(\int_0^tf(s)(|Y_s|^2+1)ds\big)^\frac{p}{2}\\&
\quad+C_pC'_p\Big(\frac{1}{2K}\mathbb{E}((Y^2_t+1)^\frac{p}{2})+\frac{K}{2}\mathbb{E}\big[\big(\int_0^tf^2(s)
(|Y_s|^2+1)ds\big)^\frac{p}{2}\big]\Big)\\&
\le C_p|x_0|^p+\frac{C_pC'_p C^{''}_p}{2K}(\mathbb{E}(Y^p_t)+1)+C_p\mathbb{E}\big(\int_0^tf(s)(|Y_s|^2+1)ds\big)^\frac{p}{2}\\&
\quad+\frac{C_pC'_pK}{2}\mathbb{E}\big[\big(\int_0^tf^2(s)
(|Y_s|^2+1)ds\big)^\frac{p}{2}\big]\Big).\endaligned\]
Taking $K:=C_pC'_pC^{''}_p,$ we have

\[\aligned\mathbb{E}(Y_t^p)&\le 1+2C_p|x_0|^p+2C_p\mathbb{E}\big(\int_0^tf(s)(|Y_s|^2+1)ds\big)^\frac{p}{2}\\&
\quad+C^2_pC'^2_pC^{''}_p\mathbb{E}\big[\big(\int_0^tf^2(s)(|Y_s|^2+1)ds\big)^\frac{p}{2}\big]\Big).\endaligned\]

If $f\in L_{loc}^p(\mathbb{R}_+),$ then

\[\aligned\mathbb{E}\Big[\big(\int_0^tf^r(s)(|Y_s|^2+1)ds\big)^\frac{p}{2}\Big]
&\le C^{''}_p\Big((\int_0^tf^r(s)ds)^\frac{p}{2}+\mathbb{E}\big[(\int_0^tf^r(s)Y_s^2ds)^\frac{p}{2}\big]\Big)\\&
\le C^{''}_p\Big((\int_0^tf^r(s)ds)^\frac{p}{2}+\int_0^tf^\frac{rp}{2}(s)\mathbb{E}(Y_s^p)ds\Big)\endaligned\]
where $r=1, 2$. The last inequality holds because of H\"{o}lder's inequality. Thus,

\[\mathbb{E}(Y_t^p)\le A+B\int_0^tf^\frac{p}{2}(s)\mathbb{E}(Y_s^p)ds+C\int_0^tf^p(s)\mathbb{E}(Y_s^p)ds\]
where

$$A:=1+2C_p|x_0|^p+2C_pC^{''}_p(\int_0^tf(s)ds)^\frac{p}{2}+C^2_p{C'}^2_p{C^{''}}^2_p(\int_0^tf^2(s)ds)^\frac{p}{2}$$
$$B:=2C_pC^{''}_p,\quad C:=C^2_p{C'}^2_p{C^{''}}^2_p.$$

By Gronwall's lemma, we have

\[\mathbb{E}(Y_t^p)\le A\exp\{B\int_0^tf(s)^\frac{p}{2}ds+C\int_0^tf^p(s)ds\}.\]

On the other hand, if $f\in L_{loc}^\frac{2p}{p-2}(\mathbb(R)_+),$ then by H\"older's inequality,

\[\aligned\mathbb{E}\Big[\big(\int_0^tf^r(s)(|Y_s|^2+1)ds\big)^\frac{p}{2}\Big]
&\le C^{''}_p\Big((\int_0^tf^r(s)ds)^\frac{p}{2}+\mathbb{E}\big[(\int_0^tf^r(s)Y_s^2ds)^\frac{p}{2}\big]\Big)\\&
\le C^{''}_p\Big((\int_0^tf^r(s)ds)^\frac{p}{2}+(\int_0^tf(s)^\frac{rp}{p-2}ds)^\frac{p-2}{p}\int_0^t\mathbb{E}(Y_s^p)ds\Big) \endaligned\]
where $r=1,2.$ So we arrive at

\begin{equation}\mathbb{E}(Y_t^p)\le A_1+B_1\int_0^t\mathbb{E}(Y_s^p)ds
\end{equation}
where

$$A_1:=A=1+2C_p|x_0|^p+2C_pC^{''}_p(\int_0^tf(s)ds)^\frac{p}{2}+C^2_p{C'}^2_p{C^{''}}^2_p(\int_0^tf^2(s)ds)^\frac{p}{2},$$
and
$$B_1:=2C_pC^{''}_p(\int_0^tf^\frac{p}{p-2}(s)ds)^\frac{p-2}{p}+C^2_p{C'}^2_p{C^{''}}^2_p(\int_0^tf^\frac{2p}{p-2}(s)ds)^\frac{p-2}{p}.$$
By Gronwall's lemma, we have

\[\mathbb{E}(Y_t^p)\le A_1e^{B_1t}.\]
$\square$

\noindent\textbf{Example}\quad We just consider the time-homogeneous case for simplicity.
Suppose $d=2,$ $m=1$. For any $r>0,$ define

$$\sigma(x)=|x|^r(-x_2,x_1)^T, b(x)=-|x|^{2r}x^T.$$

It's obvious that there exists a unique strong solution for the giving stochastic differential equation
since the local Lipschitzian condition holds for both $\sigma$ and $b.$ On the other hand,

\[\Big(|\sigma(x)|^2+2\langle x,b(x)\rangle\Big)\vee |\sigma^T(x)x|=(|x|^{2r+2}-2|x|^{2r+2})\vee 0=0\le K(|x|^2+1)\]
So by Theorem \ref{dingli4}, we can get the upper bound of the $p$th moment of the maximum process. But there is no $K>0$ such that

$$|\sigma(x)|^2=|x|^{2r+2}\le K(|x|^2+1).$$
So we have given a sufficient condition for the boundedness of the $p$th moment of the maximum process, which is weaker than that of \cite{HMS,Mao}.

\section{Non confluence of the solution}\vskip0.2in

\noindent\textbf{Proof of Theorem \ref{dingli3}}\quad Assume that
$X_t(x_0) $ is a solution of equations (\ref{sdelan1}) starting from
$x_0 $. Without loss of generality, we assume that
$0<\varepsilon<|x_0-y_0|<c_0/2 $, then define
\begin{equation}\label{tingshi4.1}\hat{\tau}_\varepsilon:=\inf\{t>0,|X_t(x_0)-X_t(y_0)|\le\varepsilon\},\
\hat{\tau}:=\inf\{t>0,X_t(x_0)=X_t(y_0)\}. \end{equation}
Clearly,
$\hat{\tau}_\varepsilon\rightarrow\hat{\tau}$, as $\varepsilon\rightarrow0.$

Let
\begin{equation}\label{tingshi4.2}\aligned&\tau:=\inf\{t>0,|X_t(x_0)-X_t(y_0)|\ge\frac{3}{4}c_0\},\\&
\tau_R:=\inf\{t>0,|X_t(x_0)|\vee|X_t(y_0)|\ge R\}\endaligned \end{equation}
and
$$Y_t:=X_{t\wedge\hat{\tau}_\varepsilon\wedge\tau_R}(x_0)-X_{t\wedge\hat{\tau}_\varepsilon\wedge\tau_R}(y_0),\
 \xi_t:=|Y_t|^2.$$

We take the test function
\begin{equation}\Phi_\delta(x):=\exp(\int_x^{c_0}\frac{ds}{\gamma_R(s)+\delta}). \end{equation}

By It\^{o}'s formula, we have
$$
\aligned\Phi_{\delta}(\xi_{t\wedge\tau})&=\Phi_{\delta}(\xi_0)+\int_0^{t\wedge\tau}\Phi'_{\delta}(\xi_s)d\xi_s
+\frac{1}{2}\int_0^{t\wedge\tau}\Phi{''}_{\delta}(\xi_s)d\langle\xi,\xi\rangle_s\\&=
\Phi_{\delta}(\xi_0)+M_t+\int_0^{t\wedge\tau}\Phi_{\delta}(\xi_s)(\frac{-1}{\gamma_R(\xi_s)+\delta})( 2\langle
Y_s,f_s\rangle
+||e_s||^2)ds\\&\quad
+\frac{1}{2}\int_0^{t\wedge\tau}\Phi_{\delta}(\xi_s)\frac{\gamma_R'(\xi_s)+1}{(\gamma_R(\xi_s)+\delta)^2}\cdot 4\xi_s||e_s||^2ds
\endaligned
$$
where
$$M_t:=2\int_0^{t\wedge\tau}\Phi'_{\delta}(\xi_s)\langle Y_s,e_sdB_s\rangle,$$
$$e_s:=\sigma(X_s(x_0))-\sigma(X_s(y_0)),$$
$$h_s:=b(X_s(x_0))-b(X_s(y_0)),$$
then taking expectation on both sides, we get
$$
\aligned\mathbb{E}(\varphi_{\delta}(\xi_{t\wedge\tau}))&=\Phi_{\delta}(\xi_0)+
\mathbb{E}\int_0^{t\wedge\tau}\Phi_{\delta}(\xi_s)\Big[\frac{2\xi_s(\gamma_R'(\xi_s)+1)||e_s||^2}{(\gamma_R(\xi_s)+\delta)^2}-\frac{2\langle
Y_s,h_s\rangle+||e_s||^2}{\gamma_R(\xi_s)+\delta}\Big]ds\\&
\le\Phi_{\delta}(\xi_0)+
\mathbb{E}\int_0^{t\wedge\tau}\Phi_{\delta}(\xi_s)\frac{2K||e_s||^2-(2\langle
Y_s,h_s\rangle+||e_s||^2)}{\gamma_R(\xi_s)+\delta}ds\\&
\le\Phi_{\delta}(\xi_0)+
(2K-1)\mathbb{E}\int_0^{t\wedge\tau}\Phi_{\delta}(\xi_s)\frac{||e_s||^2-\frac{2}{2K-1}\langle
Y_s,h_s\rangle}{\gamma_R(\xi_s)+\delta}ds\\&\le\Phi_{\delta}(\xi_0)+
(2K-1)\mathbb{E}\int_0^t\Phi_{\delta}(\xi_s)ds.
\endaligned
$$
The last inequality holds because of condition
(\ref{t4}). Then by Gronwall's lemma, we have

$$\mathbb{E}(\Phi_{\delta}(\xi_{t\wedge\tau}))\le\Phi_{\delta}(\xi_0)e^{(2K-1)t}.$$
Thus  \begin{equation}\label{qiwang4.1}
\aligned\mathbb{E}(\Phi_\delta(|X_{t\wedge\hat{\tau}_\varepsilon\wedge\tau\wedge\tau_R}(x_0)
-X_{t\wedge\hat{\tau}_\varepsilon\wedge\tau\wedge\tau_R}(y_0)|^2))\le
\Phi_{\delta}(\xi_0)e^{(2K-1)t}.
\endaligned\end{equation}

On the other hand,
\begin{eqnarray*}
&& \mathbb{E}(\Phi_\delta(|X_{t\wedge\hat{\tau}_\varepsilon\wedge\tau\wedge\tau_R}(x_0)
-X_{t\wedge\hat{\tau}_\varepsilon\wedge\tau\wedge\tau_R}(y_0)|^2))\\
&\ge& \mathbb{E}(\Phi_\delta(|X_{t\wedge\hat{\tau}_\varepsilon\wedge\tau\wedge\tau_R}(x_0)
-X_{t\wedge\hat{\tau}_\varepsilon\wedge\tau\wedge\tau_R}(y_0)|^2)
1_{\hat{\tau}_\varepsilon\le
t\wedge\tau\wedge\tau_R})\\
&=&\Phi_\delta(\varepsilon^2)P(\hat{\tau}_\varepsilon\le
t\wedge\tau\wedge\tau_R).
\end{eqnarray*}
Thus,
\begin{equation}
P(\hat{\tau}_\varepsilon\le
t\wedge\tau\wedge\tau_R)\le C_t\exp{(-\int_{\varepsilon^2}^{\xi_0}\frac{ds}{\gamma(s)+\delta})}
\end{equation}
where $C_r$ is independent of $R.$
Let $R\rightarrow\infty, \delta\rightarrow0, \varepsilon\rightarrow0 $ subsequently. We
have for any nonnegative $t, P(\hat{\tau}\le t\wedge\tau\wedge\zeta)=0 $. Let
$t\rightarrow\infty$, it follows that $P(\hat{\tau}\le \tau\wedge\zeta)=0 $. Therefore, $\xi_\cdot$ is
positive almost surely on the interval $[0,\tau].$  Now we define

\begin{equation}T_0:=0,\quad T_1:=\tau\wedge\zeta,\quad T_2:=\inf\{t>\tau,|X_t(x_0)-X_t(y_0)|\le
\frac{c_0}{2}\}\wedge\zeta\end{equation}
and generally

\begin{equation}\aligned &T_{2n}:=\inf\{t>T_{2n-1},|X_t(x_0)-X_t(y_0)|\le
\frac{c_0}{2}\}\wedge\zeta,\\&
 T_{2n+1}:=\inf\{t>T_{2n},|X_t(x_0)-X_t(y_0)|\ge \frac{3c_0}{4}\}\wedge\zeta.\endaligned\end{equation}
Similar to Fang and Zhang \cite{Fang}, it is clear that $T_n\rightarrow\zeta, a.s.$ as $n\rightarrow\infty$. By
definition $\xi_\cdot $ is positive almost surely on the interval
$[T_{2n-1},T_{2n}] $ . Since $X_t(x)$ is stochastic continuous with
 respect to the initial value $x$ (see Theorem 3, \cite{Lan}),
 according to Corollary 5.3, \cite{Fang}, the diffusion process $X_t(x)$ is Feller.
By pathwise uniqueness of solutions,
$\{X_t\}_{t\ge0}$ has the strong Markovian property (see Mao \cite{Mao}). Starting from
$T_{2n} $ and applying the same arguments as in the first part of
the proof, one can show that $\xi_\cdot $ is positive almost surely
on the interval $[T_{2n},T_{2n+1}] $, this ends the proof. $\square$

\textbf{Acknowledgement} The authors would like to thank Professor Feng-Yu Wang for his useful suggestions.


\begin{thebibliography}{99}

\bibitem{Emery} Emery, M., Non confluence des solutions dune equation
stochastique lipschitzienne, Seminaire Proba. XV. Lecture Notes in
Mathematics, vol. 850, pp587-589, Springer, Berlin Heidelberg NewYork
1981.

\bibitem{Fang} Fang, S.Z. and Zhang, T.S., A study of a class of stochastic
differential equations with non-Lipschizian coefficients, Probab.
Theory. Relat. Fields, 2005, 132:3, 356-390.

\bibitem{Fang1} Fang, S.Z. and Zhang, T.S., Stochastic differential equations
with non-Lipschitz coefficients: pathwise uniqueness and no
explosion, C.R. Acad. Sci. Paris Ser. I., 2003, 337, 737-740.

\bibitem{HMS} Higham, D.J., Mao, X. and Stuart, A.M., Strong convergence of
Euler-Maruyama methods for nonlinear stochastic differential equations, SIAM
J. Numer. Anal., 2002, 40:3, 1041-1063.

\bibitem{Hof} Hofmanov\'a, M. and Seidler, J., On weak solutions of stochastic
differential equations, Stoch. Anal. Appl., 2012, 30, 100-121.

\bibitem{Ikeda} Ikeda, I. and Watanabe, S., Stochastic differential equations
and diffusion processes, North-Holland, Amsterdam, 1981.

\bibitem{Krylov}  Krylov, N.V., On Kolmogorov's equations for finite
dimensional diffusions, Stochastic PDE's and Kolmogorov equations in
infinite dimensions (Cetraro, 1998), Lecture Notes in Math., vol.
1715, Springer, Berlin, 1999, pp. 1-63.

\bibitem{Krylov1} Krylov, N.V. and R\"ockner, M., Strong solutions of
stochastic equations with singular time dependent drift, Probab.
Theory Relat. Fields, 2005, 131, 154-196.

\bibitem{Lan} Lan, G.Q., Pathwise uniqueness and non-explosion of stochastic
differential equations with non-Lipschitzian coefficients, Acta.
Math. Sinica, Chinese series, 2009, 52:4,109-114.

\bibitem{Mao} Mao, X.R., Stochastic differential equations and applicatons, 2nd edition, Horwood, Chichester, 2007.

\bibitem{Rockner} Pr\'ev\"ot, C. and R\"ockner, M., A Concise
Course on Stochastic Partial Differential Equations, Lecture Notes
in Mathematics, 1905, Springer.

\bibitem{Renesse} v. Renesse, M. and Scheutzow, M., Existence and uniqueness of solutions of stochastic functional
differential equations, Random Oper. Stochastic Equations, 2010, 18:3, 267-284.

\bibitem{Revuz} Revuz, D. and Yor, M., Conetinuous martingales and Brownian
motion, Grund. der Math. Wissenschaften, 293, Springer-Verlag, 1991.

\bibitem{Wang} Shao, J., Wang, F.-Y. and Yuan, C., Harnack Inequalities for
 Stochastic (Functional) Differential Equations with Non-Lipschitzian
 Coefficients, arXiv:1208.5094.

\bibitem{Stroock} Stroock, D.W. and Varadhan, S.R.S, Multidimensional diffusion
processes, Springer-Verlag, 1979.

\bibitem{Watanabe} Watanabe, T. and Yamada, S., On the uniqueness of solutions of
stochastic differential equations II, J. Math. Kyoto Univ. 1971, 11, 553-563.

\bibitem{Yamada} Yamada, T. and Ogura, Y., On the strong comparison theorems
for solutions of stochastic differential equations, Z.
Wahrscheinlichkeitstheorie verw. Gebiete, 1981, 56, 3-19.

\bibitem{Yamada1} Yamada, T. and Watanabe, S., On the uniqueness of solutions of
stochastic differential equations, J. Math. Kyoto Univ. 1971, 11, 155-167.
\end{thebibliography}
\end{document}